# ON REGULARITY, TRANSITIVITY, AND ERGODIC PRINCIPLE FOR QUADRATIC STOCHASTIC VOLTERRA OPERATORS

## MANSOOR SABUROV

Department of Computational & Theoretical Sciences,
Faculty of Science, International Islamic University Malaysia,
P.O. Box, 141, 25710, Kuantan, Pahang, Malaysia
E-mail: msaburov@gmail.com

**ABSTRACT**

In this paper we showed an equivalence of notions of regularity, transitivity and Ergodic principle for quadratic stochastic Volterra operators acting on the finite dimensional simplex.

## Periodic Points Vs Regularity

The theory of linear operators has been well studied since the last century. The simplest nonlinear operator is a quadratic operator. A quadratic operator is a primary source for investigations of dynamical properties of population genetics. The problem of studying the behavior of trajectories of quadratic stochastic operators was stated in [11] and the application of quadratic stochastic operators to population genetics was studied in [9]. In [5], it was given a long self-contained exposition of the recent achievements and open problems in the theory of quadratic stochastic operators.

A quadratic stochastic operator arises in population genetics as follows: let us consider a population consisting of $m$ species. Let $x^{(0)} = \left(x_1^0, x_2^0, ..., x_m^0\right)$ be the probability distribution (where $x_i^0 = P(i)$ is the probability of $i \in I = \{1, 2, ..., m\}$) of species in the initial generation, and $P_{ij,k}$ is the probability that individuals in the $i$ th and $j$ th species interbred to produce an individual $k$, more precisely, $P_{ij,k}$ is the conditional probability $P(k|i,j)$ that $i$ th and $j$ th species interbred successfully, then they produce an individual $k$. We suppose that the population has the model of free population, i.e., there is no difference of "sex" and in any generation, the "parents" $i$ and $j$ are independent, i.e. $P(i,j) = P(i)P(j)$. Then the probability distribution $x^{(1)} = \left(x_1', x_2', ..., x_m'\right)$ (the state) of the species in the first generation can be found by the total probability

$$x_k' = \sum_{i,j=1}^{m} P(k|i,j) P(i,j) = \sum_{i,j=1}^{m} P_{ij,k} x_i^0 x_j^0$$

This means that the association $x^{(0)} \to x^{(1)}$ defines a map $V$ called the *evolution operator*. Let us provide the precise definition of the evolution quadratic stochastic operator.

Let $S^{m-1} = \left\{ x \in \mathbb{R}^{m-1} : \sum_{i=1}^{m} x_i = 1,\ x_i \geq 0,\ \forall i = \overline{1,m} \right\}$ be an $(m-1)$ dimensional simplex.

A mapping $V : S^{m-1} \to S^{m-1}$ defined as follows

$$x'_k \equiv (Vx)_k := \sum_{i,j=1}^{m} P_{ij,k} x_i x_j, \quad \forall k = \overline{1,m} \tag{1}$$

is said to be *a quadratic stochastic operator* where

$$x = (x_1, x_2, \ldots, x_m) \in S^{(m-1)},\ P_{ij,k} = P_{ji,k} \geq 0,\ \sum_{k=1}^{m} P_{ij,k} = 1\ \forall i,j = \overline{1,m}.$$

The population evolves by starting from an arbitrary state $x^{(0)}$, then passing to the state $x^{(1)} = V(x^{(0)})$ (in the next "generation"), then to the state $x^{(2)} = V(x^{(1)}) = V(V(x^{(0)})) = V^2(x^{(0)})$, and so on. Thus, states of the population described by the following discrete-time dynamical system

$$x^{(0)},\ x^{(1)} = V(x^{(0)}),\ x^{(2)} = V^2(x^{(0)}),\ x^{(3)} = V^3(x^{(0)}),\ \ldots,\ x^{(n)} = V^n(x^{(0)}) \ldots$$

The main problem of the theory of the dynamical system is to classify states $x^{(0)}$ based on the behavior of the trajectory $\{x^{(n)}\}_{n=1}^{\infty}$.

It is particularly interesting when the trajectory repeats. In this case we say that $x^{(0)}$ is *a periodic point*, i.e., $x^{(0)}$ is *a periodic point* if there is the smallest positive integer $k$ such that $x^{(k)} = x^{(0)}$. The number $k$ is called *a period* of the periodic point $x^{(0)}$. A *fixed point* is a periodic point of period-1, that is, a point $x^{(0)}$ such that $x^{(1)} = x^{(0)}$.

If the trajectory of every point converges to some fixed point (the limiting fixed point might be depended on the initial point) then a mapping is called *regular*. It is clear that if a mapping is regular then it does not have any order periodic points except fixed points. It turns out that, in one dimensional case, the converse statement holds true as well. More precisely, one of the fascinating results in one dimensional nonlinear dynamical system is that a mapping which maps a compact connected subset of the real line into itself is regular if and only if it does not have any order periodic points (see [1]). The most incredible result is that a mapping which maps a compact connected subset of the real line into itself is regular if and only if it does not have any period-2 points (see [1]).

It is natural to seek an analogy of these incredible results in the high dimensional case. However, in general, these results do not hold true in the high dimensional case. As a counter example we can consider the following quadratic stochastic operator $V : S^2 \to S^2$

$$V : \begin{cases} x'_1 = x_1^2 + 2x_1 x_2 \\ x'_2 = x_2^2 + 2x_2 x_3 \\ x'_3 = x_3^2 + 2x_1 x_3 \end{cases} \tag{2}$$

It is easy to check that this operator has fixed points $e_1 = (1,0,0)$, $e_2 = (0,1,0)$, $e_3 = (0,0,1)$, and $c = \left(\frac{1}{3}, \frac{1}{3}, \frac{1}{3}\right)$ and it does not have any order periodic points. However, the trajectory of any point of the interior $\text{int } S^2$ of the simplex except $c = \left(\frac{1}{3}, \frac{1}{3}, \frac{1}{3}\right)$ does not converge (see [3-4], [8], [9], [12]). More interestingly, the arithmetic mean (or the Cesaro mean) of the trajectory of the operator (2) does not converge (see [2], [13]). Surprisingly, any order arithmetic mean (or any order Cesaro mean) of the trajectory of the operator (2) does not converge (see [10]).

The regularity of the nonlinear operator acting on the high dimensional space could not be described in term of an absence of periodic points. Therefore, the study of the regularity is independent of interest.

## Regularity, Transitivity, and Ergodic Principle

The regularity problem was concerned for the operator (1) in [7]. Namely, it was studied the following problem: find the number $\alpha_m > 0$ such that $P_{ij,k} > \alpha_m$, $\forall i,j,k = \overline{1,m}$ implies the regularity of quadratic stochastic operators (1). It was shown [7] that if $\alpha_m = \frac{1}{2m}$ then the operator (1) is regular. The main problem is to find the smallest positive number among all $\alpha_m$ for fixed $m$ (if any) such that any quadratic stochastic operator under the condition $P_{ij,k} > \alpha_m$, $\forall i,j,k = \overline{1,m}$ is regular. One can easily check that if $m = 2$ then $\inf \alpha_2 = \frac{1}{2}(3 - \sqrt{7})$. If $m \geq 3$ the problem remains open.

However, the regularity problem was studied intensively for another class of quadratic stochastic operators [3-4] which could not be covered by previous cases.

**Definition [3-4]** An operator (1) is called *a quadratic stochastic Volterra operator* if
$$P_{ij,k} = 0, \ k \notin \{i,j\} \text{ for any } \forall i,j = \overline{1,m}.$$

Any quadratic stochastic Volterra operator can be written in the following form
$$(Vx)_k \equiv x'_k := x_k \left(1 + \sum_{i=1}^m a_{ki} x_i\right), \tag{3}$$

where $A_m = (a_{ki})_{k,i=1}^m$ is a skew-symmetric matrix with $a_{ki} \in [-1,1]$.

A nonlinear stochastic Volterra operator was studied in [6].

**Definition [3-4]** A skew-symmetric matrix $A_m$ is called *transversal* if all even order leading (principal) minors are nonzero. A quadratic stochastic Volterra operator (3) is called *transversal* if the corresponding skew-symmetric matrix $A_m$ is transversal.

**Theorem [3-4]** *The set of all transversal quadratic stochastic Volterra operators is an open, everywhere-dense subset of the set of all quadratic stochastic Volterra operators.*

In the sequel, we will only consider transversal quadratic stochastic Volterra operators without mentioning transversality.

The main approach to study the dynamics of quadratic stochastic Volterra operator is to construct its fixed points' chart by mean of tournaments [3-4].

*A tournament* is a complete directed graph. A tournament is called *transitive* if it does not contain a cycle of length 3.

It is clear that if a skew-symmetric matrix $A_m$ corresponding to the Volterra operator (3) is transversal then $a_{ij} \neq 0, \forall i \neq j$. Therefore, we can construct a tournament corresponding to a transversal Volterra operator (3) as follows [3-4]: a tournament $T_m$ consists of $m$ vertices and an edge directs form $i$ to $j$ if $a_{ij} < 0$ otherwise it directs form $j$ to $i$.

A Volterra operator (3) is called *transitive* if the corresponding tournament is transitive.

**Theorem [3-4]** *If a quadratic stochastic Volterra operator (3) is transitive then it is regular.*

In this paper we want to prove the converse statement of this theorem. Moreover, we describe the regularity of quadratic stochastic Volterra operators in term of the Ergodic principle. In order to this end we introduce the notion of the Ergodic principle for quadratic stochastic operators.

Let us fix a norm $\|x\| = \sqrt{\sum_{i=1}^{m} x_i^2}$ in the Euclidean space $\mathbb{R}^m$. For every $x \in S^{m-1}$ we define its support as follows

$$\text{supp}(x) = \{i : x_i \neq 0\}.$$

**Definition** We say that a quadratic stochastic operator $V : S^{m-1} \to S^{m-1}$ given by (1) satisfies *the Ergodic principle* on the simplex $S^{m-1}$ if for any $x, y \in S^{m-1}$ with $\text{supp}(x) = \text{supp}(y)$ one has $\lim_{n \to \infty} \|V^n x - V^n y\| = 0$, where $V^n(\cdot) = \underbrace{V(V(\ldots(V(\cdot))\ldots))}_{n}$ is $n$ times compositions of $V$.

The following theorem is the main result.

**Theorem** *Let $V: S^{m-1} \to S^{m-1}$ be a quadratic stochastic Volterra operator given by (3). The following statements are equivalent:*
  i)   *$V$ is regular;*
  ii)  *$V$ is transitive;*
  iii) *$V$ satisfies the Ergodic principle;*
  iv)  *One has $\lim_{n \to \infty} \|V^n x - V^{n+1} x\| = 0$ for any $x \in S^{m-1}$.*

## ACKNOWLEDGMENT

The author wishes to express his gratitude to Professor Rasul Ganikhodjaev for suggesting the problem and for many stimulating conversations.